\documentclass[12pt]{amsart}

\usepackage{amsmath,amscd,amssymb,latexsym, color}

\newcommand{\mat}{\begin{pmatrix}}
\newcommand{\emat}{\end{pmatrix}}

\renewcommand{\t}{\tau}

\renewcommand{\i}{\infty}
\newcommand{\G}{\Gamma}
\newcommand{\g}{\gamma}
\newcommand{\ve}{\varepsilon}

\renewcommand{\i}{\infty}
\newcommand{\lt}{\left}
\newcommand{\rt}{\right}

\newcommand{\Z}{\mathbb Z}

\newcommand{\m}{\mathbb M}
\renewcommand{\H}{\mathbb H}

\newcommand{\M}{\mathcal M}
\newcommand{\W}{\mathcal W}
\newcommand{\F}{\mathcal F}
\newtheorem{thm}{Theorem}

\theoremstyle{remark}

\numberwithin{equation}{section}
\numberwithin{thm}{section}
\begin{document}

\title[Cycle integrals of a sesqui-harmonic Maass form of weight zero]{Cycle integrals of a sesqui-harmonic Maass form of weight zero}

\author{Daeyeol Jeon, Soon-Yi Kang and Chang Heon Kim}
\address{Department of Mathematics Education, Kongju National University, Kongju, 314-701, Korea}
\email{dyjeon@kongju.ac.kr}
\address{Department of Mathematics, Kangwon National University, Chuncheon, 200-701, Korea} \email{sy2kang@kangwon.ac.kr}
\address{Department of Mathematics, Hanyang University,
 Seoul, 133-791 Korea}
\email{chhkim@hanyang.ac.kr}

\begin{abstract} Borcherds-Zagier bases of the spaces of weakly holomorphic modular forms of weights $1/2$ and $3/2$ share the Fourier coefficients which are traces of singular moduli. Recently, Duke, ~Imamo\={g}lu, and T\'{o}th have constructed a basis of the space of weight $1/2$ mock modular forms, each member in which has Zagier's generating series of traces of singular moduli as its shadow. They also showed that Fourier coefficients of their mock modular forms are sums of cycle integrals of the $j$-function which are real quadratic analogues of singular moduli. In this paper, we prove the Fourier coefficients of a basis of the space of weight $3/2$ mock modular forms are sums of cycle integrals of a sesqui-harmonic Maass form of weight zero whose image under hyperbolic Laplacian is the $j$-function. Furthermore, we express these sums as regularized inner products of weakly holomorphic modular forms of weight $1/2$.
\end{abstract}
\maketitle
\renewcommand{\thefootnote}%
             {}

\section{Introduction}
Fourier coefficients of half-integral weight modular forms carry rich information of number theoretic objects. Classical results include representation numbers of quadratic forms, partition functions and class numbers of imaginary quadratic number fields and more recent results connect the Fourier coefficients of half integral weight cusp forms with central values of quadratic twists of modular $L$-functions (see \cite{Kohnen,KZ,Wa}). The new development of the theory of more general automorphic forms has revealed that Fourier coefficients of weakly holomorphic modular forms and harmonic weak Maass forms of half-integral weights also convey various arithmetic properties (see \cite{BO1, BF,BO,DIT,DIT2,Zagier} for example and see \cite{BS} for more references). 

For instance, weakly holomorphic modular forms of weights $1/2$ and $3/2$ discussed by Borcherds \cite{B} and Zagier \cite{Zagier} are related with traces of singular moduli of the classical $j$-invariant 
$$j(\t)=q^{-1}+744+196884q+21493760q^2+\cdots,$$
where $q:=e^{2\pi i \t}=:e(\t)$ and $\t\in\H$, the upper half of the complex plane. It follows from the theory of complex multiplication that if $\t$ is an imaginary quadratic irrational number, then $j(\t)$ is an algebraic integer, called singular modulus. Let $J=j-744$, the normalized Hauptmodul for $\G=PSL_2(\Z)$ and for $k\in \mathbb Z+\frac{1}{2}$, let $M^!_k$ denote the space of weakly holomorphic modular forms of weight $k$ on $\G_0(4)$, in which each form satisfies Kohnen's plus space condition, that is, its Fourier expansion is of the form $\sum a(n)q^n $ where $a(n)$ is non-zero only for integers $n$ satisfying $(-1)^{k-1/2}n\equiv 0,1 \pmod 4$. 
Throughout, $D,d\equiv 0,1\pmod 4$. We also let $\mathcal Q_d$ denote the set of integral binary quadratic forms $Q=[a,b,c]=aX^2+bXY+cY^2$ with discriminant $d=b^2-4ac$ that are positive definite if $d<0$. For each $Q$ of negative discriminant $d$, there is a corresponding CM point $\t_Q$, the unique root of $Q(\t,1)=0$ in $\H$. As  $J(\t_Q)$ depends only on equivalence class of $Q$ under the usual linear fractional action of $\G$, we may define the twisted trace of singular moduli, for each fundamental discriminant $D>0$ and the associated genus character $\chi$, by
\begin{equation}\label{traceneg}{\rm Tr}_{d,D}(J)=\frac{1}{\sqrt{D}}\sum_{Q\in\G\backslash\mathcal Q_{dD}}\chi(Q)\frac{J(\t_Q)}{|\G_Q|},\quad (Dd<0),\end{equation}
where $\Gamma_Q$ is the group of automorphs of $Q$. 
 Zagier \cite{Zagier} showed the modularity of the generating series of traces of singular moduli by proving that for each $D>0$
\begin{equation}\label{gD} g_D(\tau)=q^{-D}-2\delta_{D, \square}-\sum_{d<0}{\rm Tr}_{d,D}(J)q^{|d|}\in M^!_{3/2},\end{equation}
where $\delta_{D, \square}=1$ if $D$ is a square and $0$ otherwise.
He also established a duality relation 
\begin{equation}\label{fd}f_d(\tau)=q^{d}+\sum_{D>0}{\rm Tr}_{d,D}(J)q^D\in M_{1/2}^!, \end{equation} for each $d<0$. Earlier in \cite{B}, Borcherds proved that $f_d(\t)\ (d<0)$ and $f_0(\t):=\theta(\t)=\sum_{n\in\Z} q^{n^2}$ have an interpretation in terms of infinite product expansions of certain meromorphic modular forms for $\G$. In fact, $\{f_d|d\leq 0\}$ and $\{g_D| D>0\}$ form bases for $M^!_{1/2}$ and $M^!_{3/2}$, respectively. 
 
Recently, Duke, Imamo\={g}lu and  T\'{o}th \cite{DIT} have extended Borcherds' basis $\{f_d|d\leq 0\}$ for $M^!_\frac12$ to a basis $\{f_d\}$ for  $\m^!_\frac12$, where $\m^!_k$ denotes the space of weight $k$ mock modular forms on $\G_0(4)$ satisfying the plus space condition.
 For each $d>0$, they constructed a unique mock modular form $f_d(\t)$ of weight $1/2$ with shadow $g_d(\t)$ having a Fourier expansion of the form  \begin{equation}\label{mockf}f_d(\t)=\sum_{D>0} a(D,d)q^D,\end{equation}
which implies that $f_d$ can be completed to a harmonic weak Maass form by addition of the non-holomorphic Eichler integral of $g_d$. Furthermore, they showed that for non-square $dD$ with both $d$ and $D$ positive, the Fourier coefficients $a(D,d)$ of $q^D$ in $f_d(\t)$ are sums of cycle integrals of $J$-function which are real analogues of traces of singular moduli:
\begin{equation}\label{mockfa}a(D,d)=\frac{1}{2\pi}\sum_{Q\in\G\backslash\mathcal Q_{dD}}\chi(Q)\int_{\G_Q\backslash S_Q}J(\t)\frac{d\t}{Q(\t,1)}:={\rm Tr}_{d,D}(J),\quad (D>0,d>0),\end{equation}
where the geodesic $S_Q$ is defined to be the oriented semi-circle $a|\t|^2+b{\rm Re}\t+c=0,$ directed counterclockwise if $a>0$ and clockwise if $a<0$.
 
In \cite{JKK}, the authors extended Zagier's basis $\{g_D|D> 0\}$ for $M^!_\frac32$ to a basis $\{g_D\}$ for  $\m^!_\frac32$ satisfying that for each $D\leq 0$, $g_D(\t)$ is a unique mock modular form of weight $3/2$ with shadow $f_D(\t)$ having a Fourier expansion of the form  
\begin{equation}\label{mockg}g_D(\t)=\sum_{d\leq 0} b(D,d)q^{|d|}.\end{equation}
The Fourier coefficients $b(D,d)$ in (\ref{mockg}) can be interpreted in terms of class numbers and modified traces of cycle integrals of a sesqui-harmonic Maass form. A sesqui-harmonic Maass form may not be annihilated by hyperbolic Laplacian, but transferred to a weakly holomorphic modular form by the operator (see Section \ref{hwmf} for a precise definition).
We first construct an infinite family of sesqui-harmonic Maass forms of weight $0$ whose images under the hyperbolic Laplacian $\Delta_0$ are the Faber polynomial $j_m$'s, which form a basis for the space of weight $0$ weakly holomorphic modular forms.
\begin{thm} \label{main2}  For each positive integer $m$, let $\hat{J}_m(\tau,s)$ be the sesqui-harmonic Maass form defined in (\ref{hatj}). If we set $\hat{J}_m(\tau):=\hat{J}_m(\tau,1)$, then we have
$$\Delta_{0}(\hat{J}_m)=-j_m -24\sigma (m),$$
where $\sigma(m)$  denotes the sum of positive divisors of $m$.\end{thm}
 
We now represent the Fourier coefficients of $g_D(\t)$ in terms of traces of cycle integrals of  $\hat{J}_1(\tau)$.
\begin{thm} \label{main1} Let $d$ and $D$ be negative discriminants. If $D$ is fundamental and $dD$ is non-square, then the Fourier coefficient $b(D,d)$ of $q^{|d|}$ in the mock modular form $g_D(\t)$ with shadow $f_D(\t)$ given in (\ref{mockg}) satisfies 
\begin{equation}\label{mockgb}b(D,d)=192\pi H(|d|)H(|D|)-8\sqrt{dD}\, {\rm Tr}_{d,D}^*\left(\hat{J}(\tau)\right).\end{equation}
Here $H(n)$ is the Hurwitz-Kronecker class number and ${\rm Tr}_{d,D}^*\left(\hat{J}(\tau)\right)$ is the modified trace defined in (\ref{deftr}) with $\hat{J}(\tau)=\hat{J}_1(\tau)$
\end{thm}
  
Zagier's result on traces of singular moduli was generalized by Bruinier and Funke \cite{BF} and Alfes and Ehlen \cite{AE}, in which the generating function for CM traces of a harmonic weak Maass form  of weight $0$ is shown to be a mock modular form of weight 3/2 whose shadow is a theta series of weight $1/2$. In \cite{DIT,BFI}, as alluded to earlier, the generating function for traces of cycle integrals of a harmonic weak Maass form of weight $0$ is proven to be a mock modular form of weight $1/2$ with shadow a weight $3/2$ weakly holomorphic modular form. The function $\hat{J}(\tau)$ in Theorem \ref{main1} is not harmonic, but a sesqui-harmonic Maass form.  Hence Theorem \ref{main1} shows that the generating function for traces of cycle integrals of a sesqui-harmonic Maass form of weight $0$ is a mock modular form of weight $3/2$ with shadow a weight $1/2$ weakly holomorphic modular form.

Furthermore, coefficients $b(D,d)$ of mock modular forms $g_D(\t)$, ($D<0$) can be expressed as regularized inner products of weakly holomorphic modular forms $f_D$ and $f_d$. Following \cite{B1}, we may define the regularized Petersson inner product of two modular forms $f$ and $g$ of weight $k$ for $\Gamma_0(4)$ with singularities only at the cusps by 
\begin{equation}\label{innprodef}
(f, g)^{reg}=\lim_{Y\to \infty} \int_{\mathcal F_4(Y)}
f(\tau)\overline{g(\tau)} y^k \frac{dx dy}{y^2},
\end{equation}
where $\mathcal F_4(Y)$ is the standard truncated fundamental domain for $\Gamma_0(4)$ obtained by removing $Y$-neighborhoods of the cusps. In \cite{DIT2}, Duke, Imamo\=glu and T\'oth showed that
for a positive fundamental discriminant $D$,
$$
(g_0, g_D)^{reg}=-\frac 34 \cdot \frac{\log \varepsilon_D}{\pi \sqrt D} \cdot h(D)
$$
where $\varepsilon_D$ is the smallest unit $> 1$ and $h(D)$ is the class number of the quadratic field $\mathbb Q (\sqrt D)$ in the narrow sense. They also showed that for positive discriminants $d$ and $D$ with non-square $dD$,
$$
(g_d, g_D)^{reg}= -\frac 34 {\rm Tr}_{d,D} (J).
$$
We establish analogous results for inner products of two different $f_d$'s.
\begin{thm} \label{main3} Let $d$ and $D$ be negative discriminants. If $D$ is fundamental and $dD$ is non-square, then we have
\begin{enumerate}
\item $\displaystyle{(f_D, f_d)^{reg}= -12\sqrt{Dd} {\rm Tr}_{d,D}^* \left(\hat{J}(\tau)\right)+288\pi H(|D|)H(|d|).}$
\item $\displaystyle{(f_0, f_d)^{reg}=-24\pi H(|d|).}$
\end{enumerate}\end{thm}

\noindent Theorem \ref{main3} implies $b(D,d)=\frac 23(f_D, f_d)^{reg}$ when $Dd>0$ is not a square. The regularized inner product $(f_0, f_d)^{reg}$ was also  discussed by Borcherds in \cite[Corollary 9.6]{B1} in a vector valued form.  
 
This paper is organized as follows. In Section 2, we present definitions of harmonic weak Maass forms, sesqui-harmonic Maass forms and mock modular forms along with examples in terms of Niebur Poincar\'e series and prove Theorem \ref{main2}. In Section 3, we use Niebur Poincar\'e series to prove Theorem \ref{main1}.  In Section 4, we compare the family of harmonic weak Maass forms of weight $2$ having Faber polynomial $j_m$'s as their shadows with a subset of the basis for the space of weight $2$ harmonic weak Maass forms found in \cite{DIT3}. Finally, in Section 5, we prove Theorem \ref{main3}.
  
\section{Maass forms and mock modular forms}\label{hwmf}
 Throughout, $\t=x+iy$ with $y>0$. Let $k\in \frac12\Z$ and $N$ be a positive integer with $4|N$ when $k$ is not an integer.  A {\it{harmonic weak Maass form}} $h$ of weight $k$ for $\G_0(N)$ is a smooth function on $\H$ which satisfies:\\ 
(i) $h|_k \g = h$ for all $\g \in \G_0(N)$, where $|_k$ is the weight $k$ slash operator,\\
(ii) $\Delta_k(h)=0$, where $\displaystyle{\Delta_k=-y^2\lt(\frac{\partial^2}{\partial x^2}+\frac{\partial^2}{\partial y^2}\rt)+iky\lt(\frac{\partial}{\partial x}+i\frac{\partial}{\partial y}\rt)}$,\\
(iii) $h$ has at most exponential growth at all cusps.\\
 Let $H_k(N)$ denote the space of weight $k$ harmonic weak Maass forms for $\Gamma_0(N)$.  There is an antilinear differential operator $\displaystyle{\xi_k:=2iy^{k}\frac{\overline{\partial }}{\partial \bar{\t}}}$ which plays important roles in the theory of harmonic weak Maass forms and more general automorphic forms. Considering  \begin{equation}\label{Dx}\Delta_k=-\xi_{2-k}\circ \xi_k,\end{equation} 
we find that if $h\in H_k(N)$, 
then $\xi_k(h)$ is a weight $2-k$ weakly holomorphic modular form. 
In general, a weight $k$ harmonic weak Maass form has a Fourier expansion at infinity of the form
$$h(\t)=\sum_{n\gg-\i}c_h^+(n)q^n+c_h^-(0)y^{1-k}+\sum_{0\neq n\ll \i}c_h^-(n)\G(1-k,-4\pi ny)q^n$$ so that
 \begin{equation}\label{sh}\xi_k(h)=(1-k)\overline{c_h^-(0)}-\sum_{0\neq n\ll \i}\overline{c_h^-(n)(-4\pi n)^{1-k}}q^{-n},\end{equation}
 where $\G(a,x)$ is the incomplete gamma function.
Following Zagier \cite{Z1}, we call the holomorphic part  $h^+(\t):=\sum_{n\gg-\i}c_h^+(n)q^n$ \textit{a mock modular form of weight $k$} and $g(\t):=\xi_k(h)$ \textit{the shadow} of the mock modular form $f$.  The non-holomorphic part $h^-:=h-h^+$ is then the Eichler integral of the weakly holomorphic modular form $g(\t)$. 
 
Now, we define another family of automorphic forms that includes harmonic weak Maass forms. A {\it{sesqui-harmonic Maass form}} $\F$ of weight $k$ for $\G_0(N)$ is a smooth function on $\H$ which satisfies (see \cite{BDR} for more information and references):\\ 
(i) $\F|_k \g = \F$ for all $\g \in \G_0(N)$, \\
(ii) $\Delta_{k,2}(\F)=0$, where $\displaystyle{\Delta_{k,2}=\xi_k\circ \Delta_k=-\xi_k\circ \xi_{2-k}\circ \xi_k}=\Delta_{2-k}\circ\xi_k$,\\
(iii) $\F$ has at most exponential growth at all cusps.

It is well-known that the space of weight $0$ weakly holomorphic modular forms on $\G$ has a unique basis $\{j_m|m\geq 0\}$ where $j_m$ is uniquely  determined by having the form $j_m=q^{-m}+O(q)$. For example, $j_0=1$ and $j_1=j-744=J$, the hauptmodul for $\G$. This basis can be extended to a basis for the space of weight $0$ harmonic weak Maass forms using  Poincar\'e series. 
If $\phi: \mathbb R^+ \to \mathbb C$ is a smooth function satisfying $\phi (y)=O_\varepsilon (y^{1+\varepsilon})$ for any $\varepsilon >0$ and $\G_\i$ is the subgroup of translations of $\G$, then the general Poincar\'e series 
\begin{equation}\label{gp}
G_m(\tau, \phi)=\sum_{\g\in \Gamma_\infty \backslash \Gamma}
e(m {\rm Re}(\g\tau)) \phi ({\rm Im}(\g\tau)),\quad (m\in\Z)
\end{equation}
is a smooth $\G$-invariant function on $\H$. 
Let
\begin{equation}\label{pms} \phi_{m,s} (y)=
\left\{\begin{array}{ll}
2\pi|m|^{\frac{1}{2}}y^{\frac{1}{2}}I_{s-\frac{1}{2}}(2\pi|m|y),
&m\neq 0,\\
y^s,
&m= 0,\end{array}\right.\end{equation}
with $I_{\nu}$ the usual $I$-Bessel function.
Then the Niebur Poincar\'e series $G_m(\tau, s):=G_m (\tau, \phi_{m,s})$ is defined for ${\rm Re} \, s > 1$ and satisfies 
\begin{equation}\label{npdelta}\Delta_0G_{m}(\t,s)=(s-s^2)G_{m}(\t,s).\end{equation}
As each $G_m(\t,s)$ when $m\neq 0$ has an analytic continuation to $Re(s) > 1/2$, we obtain an infinite family of weight $0$ harmonic weak Maass forms $G_m(\t,1)$ ($m\in \mathbb Z$).  The Fourier coefficients of $G_m(\tau, s)$ can be written in terms of $I$, $J$, $K$-Bessel functions and generalized Kloosterman sum which is defined by
\begin{equation}\label{kloos}
K_k(m,n;c):=\left\{
              \begin{array}{ll}
                \sum_{v (c)^*}e\left(\frac{m\bar{v}+nv}{c}\right), & \hbox{if $k\in\mathbb{Z}$,} \\
                \sum_{v (c)^*}\lt(\frac cv\rt)^{2k}\ve_v^{2k}e\left(\frac{m\bar{v}+nv}{c}\right), & \hbox{if $k\in\frac12 \mathbb{Z}\backslash \mathbb{Z}$,}
              \end{array}
            \right.
\end{equation}
where the sum runs through the primitive residue classes modulo $c$ and $v\bar{v}\equiv 1 \pmod c$. The function $G_0(\t,s)$ is the usual Eisenstein series whose Fourier expansion is given by
\begin{equation}\label{g0}
G_0(\t,s)=y^s+\frac{\xi(2s-1)}{\xi(2s)}y^{1-s}+\sum_{n\neq 0}\frac{2y^{1/2}\pi^s}{\G(s)\zeta(s)}|n|^{s-1/2}{\sigma_{1-2s}(|n|)}K_{s-\frac12}(2\pi|n|y)e(nx),\end{equation}
where $\xi(s)=\pi^{-s/2}\G(s/2)\zeta(s)$ and $\sigma_a(n)$ is the sum of $a$th powers of positive divisors of $n$ (\cite{DIT,IK}). For $m\neq 0$, the Fourier expansion of $G_m(\t,s)$ is given in \cite{DIT,F} by
\begin{eqnarray}\label{gm}
G_m(\t,s)&=&2\pi|m|^{1/2}y^{1/2}I_{s-1/2}(2\pi|m|y)e(mx)+\frac{4\pi|m|^{1-s}\sigma_{2s-1}(|m|)}{(2s-1)\xi(2s)}y^{1-s}\cr
&&\qquad+4\pi|m|^{1/2}y^{1/2}\sum_{n\neq 0}|n|^{1/2}c_m(n,s)K_{s-\frac12}(2\pi|n|y)e(nx),\end{eqnarray}
where $$c_m(n,s)=\sum_{c>0}c^{-1}K_0(m,n,c)
\left\{\begin{array}{cc}I_{2s-1}(4\pi\sqrt{|mn|}c^{-1}),&mn<0,\\
J_{2s-1}(4\pi\sqrt{|mn|}c^{-1}),&mn>0.\end{array}\right.$$
Following \cite{DIT}, we define for positive integer $m$
\begin{equation}\label{jms}
j_m(\t,s):=G_{-m}(\t,s)-\frac{2m^{1-s}\sigma_{2s-1}(m)}{\pi^{-(s+1/2)}\G(s+1/2)\zeta(2s-1)}G_0(\t,s)\end{equation}
and apply Weil's bound to its Fourier coefficients so that we have an analytic continuation for $j_m(\t,s)$ to $Re(s) > 3/4$. As $s\to 1$, the pole in $G_0(\t,s)$ cancels with the zero in the quotient multiplied to $G_0(\t,s)$. Moreover, it follows from (\ref{gm}) and $j_m(\t)=q^{-m}+O(q)$ that
$$j_m(\t,1)=G_{-m}(\t,1)-24\sigma(m)=j_m(\t),$$ 
because a bounded harmonic function is constant. (See \cite{F,Neun,Niebur,DIT} for more details and references for Niebur Poincar\'e series.)  

Now we construct sesqui-harmonic Maass forms $\hat{J}_m(\tau)$ satisfying $\Delta_0(\hat{J}_m(\tau))=-j_m -24\sigma (m)$.
For positive integer $m$ and $\rm{Re}\ s>1$, define  \begin{equation}\label{hatj}\hat{J}_m(\tau, s):=\frac{\partial}{\partial s}G_{-m}(\t,s)=G_{-m}(\t, \frac{\partial}{\partial s}\phi_{-m,s}).\end{equation} Then $\hat{J}_m(\tau, s)$ also has an analytic continuation to $Re(s) > 1/2$. Hence it follows from (\ref{npdelta}) that 
$$\Delta_0(\hat{J}_m(\tau))=\Delta_0\left(\frac{\partial}{\partial s}G_{-m}(\t,s)|_{s=1}\right)=-G_{-m}(\t,1)=-j_m(\t) -24\sigma (m),$$ which proves Theorem \ref{main2}.

By means of (\ref{Dx}), we observe that $h^*_m(\t):=\xi_0(\hat{J}_m(\tau))$ is a harmonic weak Maass form of weight $2$ satisfying $\xi_2(h^*_m(\tau))=j_m(\t) +24\sigma (m)$. In other words, the holomorphic part of $h^*_m(\tau)$ is a mock modular form of weight $2$ with shadow $j_m(\t) +24\sigma (m)$. Recently, Duke, Imamo\={g}lu and T\'{o}th \cite{DIT3} have found harmonic weak Maass forms $h_m(\t)$ which have bounded holomorphic parts and satisfy $\xi_2(h_m)=\frac{1}{4\pi}(j_m + 24\sigma (m))$. This implies that for each positive integer $m$, $h^*_m(\t)$ and $4\pi h_m(\t)$ differ by a weakly holomorphic modular form. In fact, $h^*_m(\t)=4\pi h_m(\t)$ and we confirm this in Section 4.

\section{Proof of Theorem \ref{main1}}\label{proofmain1}
For negative $dD$, ${\rm Tr}_{d,D}(J)$ is a Fourier coefficient of a weakly holomorphic modular form while, for positive discriminants $D$ and $d$ with $dD$ non-square, ${\rm Tr}_{d,D}(J)$ is a Fourier coefficient of a harmonic weak Maass form. But when both $D$ and $d$ are negative discriminants,  ${\rm Tr}_{d,D}(J)$ is simply zero, as seen below.  

Suppose that $D$ and $d$ are discriminants and $D$ is fundamental. For a quadratic form $Q=[a,b,c]$ with discriminant $dD$, define a genus character $\chi$ on $\Gamma\backslash Q_{dD}$ by (\cite{GKZ})
$$\chi_D(Q)=\left\{\begin{array}{cc}(\frac{D}{r}),& (a,b,c,D)=1\ \mathrm{and}\ (r,D)=1\ \mathrm{where}\ Q\ \mathrm{represents}\ r,\\
0, &  (a,b,c,D)>1.\end{array}\right.$$ 
Then for a non-square $dD$, the twisted trace of a general Poincar\'e series defined in (\ref{gp}) is given by 
$${\rm Tr}_{d,D}(G_m(\tau, \phi))=\frac{1}{2\pi}\sum_{Q\in\G\backslash\mathcal Q_{dD}}\chi(Q)\int_{\G_Q\backslash S_Q}G_m(\tau, \phi)\frac{d\t}{Q(\t,1)}.$$
If we denote the subset of $\mathcal Q_{dD}$ consisting of quadratic forms with $a>0$ by $\mathcal Q_{dD}^+$, then it follows from the proof of \cite[Lemma 7]{DIT} that
\begin{eqnarray*}{\rm Tr}_{d,D}(G_m(\tau, \phi))&=& \sum_{\Gamma_\infty\backslash \mathcal Q_{dD}^+}
\chi (Q) \int_{S_Q} e(m {\rm Re} \, \tau) \phi ({\rm Im} \, \tau) \frac{d\tau}{Q(\tau,1)}  \cr
&+& \sum_{\Gamma_\infty\backslash \mathcal Q_{dD}^+}
\chi (-Q) \int_{S_{-Q}} e(m {\rm Re} \, \tau) \phi ({\rm Im} \, \tau) \frac{d\tau}{-Q(\tau,1)}.\cr
\end{eqnarray*}
Since the genus character $\chi$ satisfies $\chi(-Q)=\chi_D(-Q)={\rm sgn}( D) \chi_D (Q)$, we have 
 $${\rm Tr}_{d,D}(G_m(\tau, \phi))=\left\{
  \begin{array}{ll}
    2\sum_{\Gamma_\infty\backslash \mathcal Q_{dD}^+}
\chi (Q) \int_{S_Q} e(m {\rm Re} \, \tau) \phi ({\rm Im} \, \tau) \frac{d\tau}{Q(\tau,1)}, & \hbox{if $d>0, D>0$, } \\
    0, & \hbox{if $d<0, D<0$.}
  \end{array}
\right.$$ 
Thus for negative discriminants $D$ and $d$,
${\rm Tr}_{d,D}(G_m(\tau, \phi))=0,$
and hence ${\rm Tr}_{d,D}(J(\t))={\rm Tr}_{d,D}(G_{-1}(\tau,1)-24)=0.$

It is then reasonable to define a modified trace
${\rm Tr}_{d,D}^*$ for each Poincar\'e series $G_m(\tau, \phi)$ by
\begin{equation}\label{modtrps}
{\rm Tr}_{d,D}^*(G_m(\tau, \phi))=
2\sum_{\Gamma_\infty\backslash \mathcal Q_{dD}^+}
\chi (Q) \int_{S_Q} e(m {\rm Re} \, \tau) \phi ({\rm Im} \, \tau) \frac{d\tau}{Q(\tau,1)}
\end{equation}
so that ${\rm Tr}_{d,D}^*(G_m(\tau, \phi))={\rm Tr}_{d,D}(G_m(\tau, \phi))$ when both $D$ and $d$ are positive. From now to the end of this section, we assume both $d$ and $D$ are negative discriminants and $D$ is fundamental with  $dD$ non-square. Modifying
\cite[Lemma 7]{DIT}, we obtain that
$$
2\pi \sqrt{dD} {\rm Tr}_{d,D}^*(G_m(\tau, \phi))
=\sum_{0<c\equiv0\, (4)}S_m(d,D;c)\Phi_m\lt(\frac{2\sqrt{dD}}{c}\rt),
$$
where
$\Phi_m(t)=\int_0^\pi\cos(2\pi mt\cos\theta)\phi(t\sin\theta)\frac{d\theta}{\sin\theta}$ for $t>0$ and 
$$S_m(d,D;c)=\sum_{b\ (\rm{mod}\, c)\atop b^2\equiv Dd\ (\rm{mod}\,  c)}\chi\left(\lt[\frac c4,b,\frac{b^2-Dd}{c}\rt]\right)e\lt(\frac{2mb}{c}\rt).$$
If we take $\phi=\phi_{m,s}$ defined in (\ref{pms}) for nonzero integer $m$, we then have
\begin{align*}
\Phi_m(t) &=\int_0^\pi\cos(2\pi mt\cos\theta)2\pi\sqrt{|m|}(t\sin\theta)^{\frac{1}{2}}I_{s-\frac{1}{2}}(2\pi|m|t\sin\theta)
\frac{d\theta}{\sin\theta} \\
&=2\pi\sqrt{|m|t}\int_0^\pi\cos(2\pi mt\cos\theta)I_{s-\frac{1}{2}}(2\pi|m|t\sin\theta)
\frac{d\theta}{(\sin\theta)^{\frac{1}{2}}}\\
&=\pi|m|^{\frac{1}{2}}t^{\frac{1}{2}}\frac{2^s\Gamma(\frac{s}{2})^2}{\Gamma(s)} J_{s-\frac{1}{2}}(2\pi|m|t),
\end{align*}
where the last equality follows from \cite[Lemma 9]{DIT} for $Re(s)>0$.
Applying \cite[Proposition 3]{DIT} in the second equality below with $K^+(m,n;c)=(1-i)\lt(1+(\frac{4}{c/4})\rt)K_{1/2}(m,n;c)$ and making a suitable change of variables in the last, we find that
\begin{eqnarray}\label{modtra} 
&&\frac{\Gamma(s)}{2^s\Gamma(\frac s2)^{2}}2\pi \sqrt{dD} \, {\rm Tr}_{d,D}^*(G_{-m}(\tau, s))\\
&=& \pi\sqrt{2|m|}(dD)^{\frac{1}{4}}\sum_{0<c\equiv0\, (4)}\frac{S_{-m}(d,D;c)}{\sqrt{c}}J_{s-\frac{1}{2}}\left(\frac{4\pi|m|\sqrt{dD}}{c}\right)\cr
&=&\pi\sqrt{2|m|}(dD)^{\frac{1}{4}}\sum_{0<c\equiv0\, (4)}\frac{1}{\sqrt{c}}\sum_{n|(m,\frac{c}{4})}\left(\frac{D}{n}\right)
\sqrt{\frac{n}{c}}K^+\left(d,\frac{m^2D}{n^2};\frac{c}{n}\right)J_{s-\frac{1}{2}}\left(\frac{4\pi|m|\sqrt{dD}}{c}\right)\cr
&=&\pi\sqrt{2|m|}(dD)^{\frac{1}{4}}\sum_{n|m}\left(\frac{D}{n}\right)n^{-\frac{1}{2}}\sum_{0<c\equiv0\, (4)}\frac{1}{c}
K^+\left(d,\frac{m^2D}{n^2};c\right)J_{s-\frac{1}{2}}\left(\frac{4\pi}{c}\sqrt{\frac{m^2}{n^2}Dd}\right).\cr
\end{eqnarray}
On the other hand, in \cite{JKK}, the authors constructed general Maass-Poincar\'{e} series using spherical Whittaker functions $\M_{n}(y,s)$ and $\W_{n}(y,s)$ which are defined by
\begin{eqnarray*}\label{defmnwn}
&&\mathcal{M}_{n}(y,s)=\left\{
                     \begin{array}{ll}
                       \G(2s)^{-1}(4\pi |n|y)^{-k/2}M_{\frac k2\mathrm{sgn}(n),s-1/2}(4\pi |n|y), & \hbox{if $n\neq 0$,} \\
                       y^{s-k/2}, & \hbox{if $n=0$,}
                     \end{array}
                   \right.\\ \nonumber
&&\mathcal{W}_{n}(y,s)=\left\{
  \begin{array}{ll}
   \G(s+\frac k2\mathrm{sgn}(n))^{-1}|n|^{k/2-1}(4\pi y)^{-k/2}W_{\frac k2\mathrm{sgn}(n),s-1/2}(4\pi |n| y), &  \hbox{if $n\neq 0$,}\\
   \frac{(4\pi)^{1-k}y^{1-s-k/2}}{(2s-1)\G(s-k/2)\G(s+k/2)}, &  \hbox{if $n=0$,}
  \end{array}
\right.\end{eqnarray*}
where $M_{\mu,\nu}(y)$ and $W_{\mu,\nu}(y)$ are Whittaker functions. (See \cite[Section 2]{JKK} or (\ref{conf}) for definition.) 
In particular, families of Maass-Poincar\'{e} series of weight $3/2$ on $\Gamma_0(4)$ satisfying the plus space condition are given by
\begin{equation}\label{for32}F_{m}^+(\t,s)=\M_{m}(y,s)e(mx)+\sum_{n\equiv 0,3\, (4)}b_{m}(n,s)\W_{n}(y,s)e(nx)\end{equation}
for each $m\equiv 0,3 \pmod 4$ in \cite[Theorem 4.4]{JKK}. 
In case both $m$ and $n$ are positive, it follows from \cite[Theorem 4.4]{JKK} and the property of Kloosterman sum
$$K_{\frac{3}{2}}(m,n;c)=-iK_{\frac{1}{2}}(-m,-n;c)$$
that
\begin{equation}\label{bmns}b_{m}(n,s)=-\sqrt{2}\pi\sum_{0<c\equiv0\, (4)}\frac{K^+(-m,-n;c)}{c}|mn|^{-\frac{1}{4}}J_{2s-1}\left(\frac{4\pi\sqrt{|mn|}}{c}\right).\end{equation}
Comparing (\ref{bmns}) with (\ref{modtra}), we find that
\begin{equation}\label{modtrzero}\frac{\Gamma(s)}{2^s\Gamma(\frac s2)^{2}}2\pi  \, {\rm Tr}_{d,D}^*(G_{-m}(\tau, s))\\=-\sum_{n|m}\sqrt{\frac{m^2}{n^2}}\left(\frac{D}{n}\right)b_{|D|}\left(\frac{m^2}{n^2}|d|,\frac{s}{2}+\frac{1}{4}\right).\end{equation}
This implies that  ${\rm Tr}_{d,D}^*(G_{-m}(\tau, s))$ should vanish at $s=1$ as $F_m^+(\t,3/4)=\{0\}$ for positive integer $m$, according to \cite[Proposition 5.1]{JKK}. 

Now we differentiate both sides of (\ref{modtrzero}) with respect to $s$ at $s=1$ so that we have 
\begin{align*}
{\rm Tr}_{d,D}^*(\hat{J}_m(\tau,s))|_{s=1}
&=-\sum_{n|m}\left|\frac{m}{n}\right|\left(\frac{D}{n}\right)\frac{\partial}{\partial s}\lt[
b_{|D|}\left(\frac{m^2}{n^2}|d|,\frac{s}{2}+\frac{1}{4}\right)\rt]_{s=1}
\\
&=-\frac{1}{2}\sum_{n|m}\left|\frac{m}{n}\right|\left(\frac{D}{n}\right)\frac{\partial}{\partial s}\lt[b_{|D|}\left(\frac{m^2}{n^2}|d|,s\right)\rt]_{s=\frac{3}{4}}.
\end{align*}
For simplicity, if we define
\begin{equation} \label{deftr}
{\rm Tr}_{d,D}^*(\hat{J}_m(\tau))={\rm Tr}_{d,D}^*(\hat{J}_m(\tau,s))|_{s=1},
\end{equation}
then we have
\begin{equation}\label{m1}
{\rm Tr}_{d,D}^*(\hat{J}(\tau))={\rm Tr}_{d,D}^*(\hat{J}_1(\tau))
 =-\frac{1}{2}\frac{\partial}{\partial s}\lt[b_{|D|}(|d|,s)\rt]_{s=\frac{3}{4}}.\end{equation}

For each negative discriminant $D<0$, the mock modular form $g_D(\t)$ in (\ref{mockg}) and \cite[Theorem 1.1]{JKK} is given by $g_D(\t)=2\sqrt{\pi |D|}k_D^+(\t)$, where $k_D^+(\t)$ denotes the holomorphic part of a harmonic weak Maass form $h_{-D,3/2}(\t)$ found in \cite[Theorem 5.3]{JKK}.  More precisely, the harmonic weak Maass form of weight $3/2$ satisfying the plus space condition was constructed via
\begin{equation}\label{kD}k_D(\t):=h_{-D,3/2}(\t)=\frac{\partial}{\partial s}F_{|D|}^+(\t,s)|_{s=\frac{3}{4}}-8\sqrt{\frac{\pi}{|D|}}H(|D|)F_0^+(\t,3/4)\end{equation}
and it follows from \cite[Theorem 5.2 and Theorem 5.3]{JKK} that the Fourier expansion of $k_D^+$ is given by
\begin{eqnarray}\label{kplus}
k_D^+&=& -2\sqrt{\pi}i q^{|D|} -8 \sqrt{\frac{\pi}{|D|}} H(|D|)\cr &&+\sum_{0< n\equiv 0,3\, (4)\atop n\neq |D|} 
\left(
\frac{\partial}{\partial s}\lt[b_{|D|}(n,s)\rt]_{s=\frac{3}{4}} \frac{2\sqrt n}{\sqrt \pi}+96\sqrt{\frac{\pi}{|D|}}H(|D|) H(n)
\right)q^n.
\end{eqnarray}
Theorem \ref{main1} then follows from (\ref{m1}) and (\ref{kplus}).

\section{Weight $2$ harmonic weak Maass form}

In Section 2, for positive integers $m$, we established weight $0$ sesqui-harmonic Maass forms $\hat{J}_m(\tau)$  and made an observation that $h_m^*(\t)=\xi_0(\hat{J}_m(\tau))$ are harmonic weak Maass forms of weight $2$ satisfying $\xi_2(h_m^*(\tau))=j_m(\t) +24\sigma (m)$. The members $h_m$ ($m>0$) of the basis for the space of weight $2$ harmonic weak Maass forms found by Duke, Imamo\={g}lu and T\'{o}th in \cite{DIT3} have the same property. In this section, we prove $h_m$ and $h_m^*$ are equal up to a constant multiple.

By its definition, $h_m^*(\tau)=\xi_0(\hat{J}_m(\tau))$ and  
\begin{eqnarray*} \xi_0(\hat{J}_m(\tau,s))&=&\sum_{\gamma\in\Gamma_\infty\backslash\Gamma}\left[\xi_0\left(\frac{\partial}{\partial s}\phi_{-m,s}(y)e(-mx)\right)\right]\mid_2\gamma\\
&=&\sum_{\gamma\in\Gamma_\infty\backslash\Gamma}\left(\frac{\partial}{\partial s}\xi_0\left(\phi_{-m,s}(y)e(-mx)\right)\right)\mid_2\gamma.\end{eqnarray*}
It follows from \cite[(13.6.3)]{AS}, \cite[p.10]{BO3} or \cite[Appendix A]{DIT} that
$$\phi_{-m,s}(y)=2\pi\sqrt{m}\sqrt{y}I_{s-\frac{1}{2}}(2\pi my)=2^{1-2s}\Gamma(s+\frac{1}{2})^{-1}\sqrt{\pi}M_{0,s-\frac{1}{2}}(4\pi my),$$ where the $M_{\mu,\nu}(y)$ is the $M$-Whittaker function given by $$M_{\mu,\nu}(y)=e^{-y/2}y^{\nu+1/2}M(\nu-\mu+\frac12,1+2\nu,y)$$
and \begin{equation}\label{conf}
M(a,b;x)=\sum_{n=0}^\i\frac{(a)(a+1)(a+2)\cdots(a+n-1)}{(b)(b+1)(b+2)\cdots(b+n-1)}\frac{x^n}{n!}.
\end{equation}
If we set $A(s)=2^{1-2s}\Gamma(s+\frac{1}{2})^{-1}\sqrt{\pi}$, then we have
\begin{eqnarray*}\xi_0(\phi_{-m,s}(y)e(-mx))&=&\overline{A(s)}\xi_0(M_{0,s-\frac{1}{2}}(4\pi my)e(-mx))\\
&=&\overline{A(s)}\xi_0(M_{0,s-\frac{1}{2}}(4\pi my)e^{-2\pi my})\bar{q}^{-m}\\
&=&\overline{A(s)}\frac{\partial}{\partial y}(M_{0,s-\frac{1}{2}}(4\pi my)e^{-2\pi my})\bar{q}^{-m} \\
&=&\overline{A(s)}(4\pi m)
\lt(sY^{-1}M_{0,s-\frac{1}{2}}(Y)e^{-\frac{Y}{2}}-\frac{1}{2}\frac{1}{\sqrt{Y}}M_{\frac{1}{2},s}(Y)e^{-\frac{Y}{2}}\rt)\bar{q}^{-m} \\
&=&\overline{A(s)}(4\pi m)
(se^{-Y}Y^{s-1}M(s,2s;Y)-\frac{1}{2}e^{-Y}Y^{s}M(s,2s+1;Y)
)
\bar{q}^{-m}
\end{eqnarray*}
where the penultimate equality follows from \cite[p.127]{BDR} with $Y=4\pi my$ and the last equality follows from (\ref{conf}). Using the transformation formula for confluent hypergeometric function  $\displaystyle{
M(\alpha,\gamma;Y)=M(\alpha+1,\gamma;Y)-\frac{Y}{\gamma}M(\alpha+1,\gamma+1;Y)}$ from \cite[(9.9.12)]{Lebedev} and (\ref{conf}), we obtain
$$\xi_0(\phi_{-m,s}(y)e(-mx))=\overline{A(s)}(4\pi m)\frac{s}{Y}M_{1,s-\frac{1}{2}}(Y)e(mx).$$ Now, let $\varphi_{m,s}(y):=(4\pi y)^{-1}M_{1,s-\frac{1}{2}}(4\pi my)$.
Then 
\begin{equation}\label{hm}\xi_0(\phi_{-m,s}(y)e(-mx))=(\overline{A(s)}4\pi s) \varphi_{m,s}(y)e(mx).\end{equation}
Differentiating both sides with respect to $s$ and summing over all $\gamma\in\Gamma_\infty\backslash\Gamma$ after applying the weight $2$ slash operator, we derive 
$$h_m^*(\t)=4\pi\sum_{\gamma\in\Gamma_\infty\backslash\Gamma}\left(\frac{\partial}{\partial s}\left(\varphi_{m,s}(y)e(mx)\right)\right)\mid_2\gamma,$$ because $\varphi_{m,1}(y)e(mx)=0$. As
$\displaystyle{\frac{\partial}{\partial s}\varphi_{m,s}(y)}$ is  the spherical function  used to construct $h_m$ in \cite{DIT3},
$h_m^*(\t) =4\pi h_m(\t).$

\section{Proof of Theorem \ref{main3}}
\par
For a smooth function $f$ of weight $k$ for $\Gamma_0(4)$ with the Fourier expansion 
$$ f(\t)=\sum_n a(n,y) e(nx),$$ we set 
$$f^e(\t)=\sum_{n\equiv 0\, (2)} a\left(n, \frac y4\right)e\left(\frac{nx}{4}\right) \hbox{ \, \, and \, \, } f^o(\t)=\sum_{n\equiv 1\, (2)} a\left(n,\frac y4\right)e\left(\frac n8\right)e\left(\frac{nx}{4}\right). 
$$
For negative discriminants $D$ and $d$, the function $f_D$ in (\ref{fd}) is holomorphic of weight $1/2$ and $k_d$ in (\ref{kD}) is smooth of weight $3/2$ for $\Gamma_0(4)$ on $\mathbb H$, both of which have Fourier expansions satisfying the plus space condition. Let $\mathcal F_4(Y)$ be the truncated domain used in \cite{DIT2}, namely, the domain obtained from the fundamental domain for $\G_0(4)$ given in \cite[Figure 1]{DIT2} by truncating at cusp $i\i$ by the line $\textrm{Im}(\t)=Y$, at the cusp $1/2$ by the circle $|\t-(\frac12+\frac{i}{8Y})|=\frac{1}{8Y}$, and at cusp $0$ by the circle  $|\t-\frac{i}{8Y}|=\frac{1}{8Y}$. Then \cite[Lemma 2]{DIT2} implies that for $Y\ge 2$, 

\begin{eqnarray*}\label{regularization}
(f_D, \xi_{3/2} k_d)^{reg}&=& \lim_{Y\to\infty}\int_{\mathcal F_4(Y)}
f_D(\tau)\overline{\xi_{3/2}( k_d(\tau))}y^{1/2}\frac{dxdy}{y^2}\\
&=&\lim_{Y\to\infty}\int_{-\frac{1}{2}+iY}^{\frac{1}{2}+iY}\left(f_D(\tau)k_d(\tau)+\frac{1}{2}f_D^e(\tau)k_d^e(\tau)
+\frac{1}{2}f_D^o(\tau)k_d^o(\tau)\right)d\tau.
\end{eqnarray*}
Letting $\t_Y=x+iY$, we have 
\begin{eqnarray*}
(f_D, \xi_{3/2} k_d)^{reg}&=&\lim_{Y\to\infty}\int_{-\frac{1}{2}}^{\frac{1}{2}}\left(f_D(\t_Y)k_d^{-}(\t_Y)+\frac{1}{2}f_D^e(\t_Y)(k_d^{-})^e(\t_Y)
+\frac{1}{2}f_D^o(\t_Y)(k_d^{-})^o(\t_Y)\right)dx \\
&&+\lim_{Y\to\infty}\int_{-\frac{1}{2}}^{\frac{1}{2}}\left(f_D(\t_Y)k_d^{+}(\t_Y)+\frac{1}{2}f_D^e(\t_Y)(k_d^{+})^e(\t_Y)
+\frac{1}{2}f_D^o(\t_Y)(k_d^{+})^o(\t_Y)\right)dx.
\end{eqnarray*}
According to \cite[Proposition 5.2, Theorem 5.3 and Eq.~(5.12)]{JKK}, the Fourier expansion of the non-holomorphic part $k_d^-(\t)$ of $k_d(\t)$ is given by
\begin{align*}
k_d^{-}(\t)& =(-i)\Gamma(-\frac{1}{2},4\pi dy)q^{-d}+\sum_{{n<0}\atop{n\equiv 0,3\, (4)}}b_{-d}\lt(n,\frac{3}{4}\rt)\sqrt{|n|}\Gamma(-\frac{1}{2},4\pi|n|y)q^n \\ & \qquad \qquad +\frac{24H(|d|)}{\sqrt d}\sum_{{0<n=\square}}\sqrt{|n|}\Gamma(-\frac{1}{2},4\pi|n|y)q^{-n},
\end{align*}
where $b_{D}(n,s)$ is given by (\ref{bmns}). Since $\Gamma(-\frac 12,4\pi|n|y)\sim e^{-4\pi|n|y}(4\pi|n|y)^{-\frac 32}$ as $y\to \i$,
\begin{align*}
(f_D, \xi_{3/2} k_d)^{reg}=&\lim_{Y\to\infty}\int_{-\frac{1}{2}}^{\frac{1}{2}}\left(f_D(\t_Y)k_d^{+}(\t_Y)+\frac{1}{2}f_D^e(\t_Y)(k_d^{+})^e(\t_Y)
+\frac{1}{2}f_D^o(\t_Y)(k_d^{+})^o(\t_Y)\right)dx\\
=&{\rm Constant\ term\ of\ }\lt(f_D(\t)k_d^{+}(\t)+f_D^e(\t)(k_d^{+})^e(\t)+f_D^o(\t)(k_d^{+})^o(\t)\rt).\end{align*}
Hence it follows from $f_D(\t)=q^D+O(1)$ that
\begin{center}
$(f_D, \xi_{3/2} k_d)^{reg}=\frac{3}{2}\times$ coefficient of $q^{|D|}$ in $k_d^{+}(\t).$
\end{center}
Recall from Section 3 that the mock modular form $g_d(\t)=2\sqrt{\pi|d|}k_d^+(\t)$ has shadow $f_d$. Now, Theorem \ref{main3} (i) follows from Theorem \ref{main1}.
Using similar arguments and (\ref{kplus}), we obtain
\begin{align*}
(f_0,f_{d})^{reg}
&=(f_0, \xi_{\frac 32} 
(2\sqrt{\pi |d|} k_d))^{reg}
=\hbox{$\frac{3}{2}\times $ constant term in 
$2\sqrt{\pi |d|} \, k_d^{+}(\t)$}
\\
&=-24\pi H(|d|),
\end{align*} 
which proves Theorem \ref{main3} (ii).

Alternatively, we have 
\begin{align*}
(f_0,f_{d})^{reg}=
(f_d,f_{0})^{reg}
&=(f_d, \xi_{\frac 32} 
(-16\pi E))^{reg}
=\hbox{$\frac{3}{2}\times (-16\pi)\times $ coefficient of $q^{|D|}$ in 
$ E^{+}(\t)$}
\\
&=-24\pi H(|d|).
\end{align*} 
Here $E(\tau):=-\frac{1}{12}k_0(\tau)$ denotes Zagier's Eisenstein series of weight $3/2$.
\section* {Acknowledgements}


\begin{thebibliography}{99}
\bibitem{AS}
M.~Abramowitz and I.~A.~Stegun, \emph{Handbook of mathematical functions with formulas, graphs, and mathematical tables},
Dover, 1972
%
\bibitem{AE}
C.~Alfes and S.~Ehlen, \emph{Twisted traces of CM values of weak Maass forms}, preprint.
%
\bibitem{B}
R.~E.~Borcherds,
\emph{Automorphic forms on $O_{s+2,2}(R)$ and infinite products}, Invent. Math. \textbf{120} (1995), no.1, 161--213.
%
\bibitem{B1}
R.~E.~Borcherds,
\emph{Automorphic forms with singularities on Grassmannians}, Invent. Math. \textbf{132} (1998), no.3, 491--562.
%
\bibitem{BDR}
K.~Bringmann, N.~Diamantis and M.~Raum,
\emph{Mock period functions, sesquiharmonic Maass forms, and non-critical values of L-functions}, Advances in Math. \textbf{233} (2013), 115--134.
%
\bibitem{BO1}
K.~Bringmann and K.~Ono,
\emph{Arithmetic properties of coefficients of half-integral weight Maass-Poincar\'e series}, Math.~Ann.~
\textbf{337} (2007), 591--612.
%
\bibitem{BF}
J.~H.~Bruinier and J.~Funke, \emph{Traces of CM-values of modular
functions}, J.~Reine Angew.~Math.~\textbf{594} (2006), 1--33.
%
\bibitem{BFI}
J.~H.~Bruinier, J.~Funke and \"{O}.~Imamo\={g}lu, \emph{Regularized theta liftings and periods of modular functions}, preprint.
%
\bibitem{BO}
J.~H.~Bruinier and K.~Ono, \emph{Heegner divisors, L-functions and harmonic weak Maass forms}, Annals of Math. \textbf{172} (2010), 2135--2181.
%
\bibitem{BO3}
J.~H.~Bruinier and K.~Ono,
\emph{Algebraic formulas for the coefficients of half integral weight harmonic weak Maass forms}, Preprint.
%
\bibitem{BS}
J.~H.~Bruinier and F.~Str\"omberg, \emph{Computation of harmonic weak Maass forms},  Experimental Mathematics, 
\textbf{21} (2012), no.~2, 117--131.
%
\bibitem{DIT}
W.~Duke, \"{O}.~Imamo\={g}lu, and \'{A}.~T\'{o}th, \emph{Cycle integrals of the $J$-function and mock modular forms},  Ann. Math. \textbf{173} (2011), no.~2, 947--981
%
\bibitem{DIT2}
W.~Duke, \"{O}.~Imamo\={g}lu, and \'{A}.~T\'{o}th, \emph{Real quadratic analogues of traces of singular invariants},  IMRN 2011, no.~13, 3082--3094
%
\bibitem{DIT3}
W.~Duke, \"{O}.~Imamo\={g}lu, and \'{A}.~T\'{o}th, \emph{Weight two weakly harmonic Maass forms}, Preprint.
%
\bibitem{F}
J.~Fay,
\emph{Fourier coefficients of the resolvent for a Fuchsian group}, J.~reine angew.~Math. \textbf{294} (1977), 143--203.
%
\bibitem{GKZ}
{B.~Gross, W.~Kohnen, and D.~Zagier}, \emph{Heegner points and derivatives of $L$-series, II}, Math. Ann. \textbf{278} (1987), 497--562.
%
\bibitem{IK}
H.~Iwaniec and E.~Kowalski, \emph{Analytic Number Theory},  American Mathematical Society Colloquium Publications \textbf{53}, American Mathematical Society, Providence, R.~I., 2004.
%
\bibitem{JKK}
D.~Jeon, S.-Y.~Kang, and C.H.~Kim, \emph{Weak Maass-Poincar\'{e} series and weight $3/2$ mock modular forms}, Preprint.
%
\bibitem{Kohnen}
W.~Kohnen, \emph{Fourier coefficients of modular forms of half-integral weight}, Math. Ann. \textbf{271} (1985), 237--268.
%
\bibitem{KZ}
W.~Kohnen and D.~Zagier, \emph{Values of L-series of modular forms at the center of the critical strip}, Invent.~Math. \textbf{64} (1981), no.~2, 175--198.
%
\bibitem{Lebedev}
N.~Lebedev, \emph{Special functions and their applications},
Dover, 1972.
%
\bibitem{Niebur}
D.~Niebur, \emph{A class of nonanalytic automorphic functions}, Nagoya Math.~J. \textbf{52} (1973), 133--145.
%
\bibitem{Neun}
H.~Neunh\"{o}ffer, \emph{\"{U}ber die analytische Fortsetzung von Poincarreihen}, (German) S.-B. Heidelberger Akad.~Wiss.~Math.-Natur.~Kl. 1973, 33--90.
%
\bibitem{Wa} 
J.-L.~Waldspurger, \emph{Sur les coefficients de Fourier des formes modulaires de poids demi-entier}, J.~Math.~Pures Appl. \textbf{60} (1981),  no.~4, 375--484.
%
\bibitem{Zagier}
D.~Zagier, \emph{Traces of singular moduli}, Motives, polylogarithms and Hodge theory, Part I, Irvine, CA,
1998, International Press Lecture Series 3, Part I (Int. Press, Somerville, MA, 2002) 211–-244.
%
\bibitem{Z1}
D.~Zagier, \emph{Ramanujan's mock theta functions and their applications},  S\'{e}minaire Bourbaki 60\`{e}me ann\'{e}e (986) (2006-2007), http://www.bourbaki.ens.fr/TEXTES/986.pdf.
%

\end{thebibliography}
\end{document}